\documentclass[11pt]{article}
\usepackage{latexsym}
\usepackage{amsfonts}
\usepackage{amsmath}
\usepackage{epsfig}

\newtheorem{thm}{Theorem} \newtheorem{lemma}[thm]{Lemma}
\newtheorem{prop}[thm]{Proposition} 

\newcommand{\prob}{\mbox{\bf P}}
\newcommand{\E}{{\bf E\,}}

\newcommand{\eps}{\epsilon}
\newcommand{\be}{\begin{equation}}
\newcommand{\ee}{\end{equation}}
\newcommand{\C}{{\cal{C}}}
\newcommand{\F}{{\cal{F}}}

\def\qed{\relax\ifmmode\hskip2em \Box\else\unskip\nobreak\hfill
$\Box$\fi}

\newcounter{mycount}
\newenvironment{proof}{\noindent{\sc Proof. }}{{\hfill
$\Box$}\par\vskip2\parsep}

\title{Component sizes of the random graph outside the scaling window}
\author{{\sc Asaf Nachmias and Yuval Peres\thanks{Microsoft Research and U.C. Berkeley. Research of both authors supported in part by NSF grants \#DMS-0244479 and \#DMS-0104073}}}

\begin{document}
\maketitle

\begin{abstract}
We provide simple proofs describing the behavior of the largest
component of the Erd\H{o}s-R\'enyi random graph $G(n,p)$ outside
of the scaling window, $p={1+\eps(n) \over n}$ where $\eps(n) \to
0$ but $\eps(n)n^{1/3} \to \infty$.
\end{abstract}

\section{Introduction}

Consider the random graph $G(n,p)$ obtained from the complete
graph on $n$ vertices by retaining each edge with probability $p$
and deleting each edge with probability $1-p$. We denote by $\C
_j$ the $j$-th largest component. Let $\eps(n)$ be a non-negative
sequence such that $\eps(n) \to 0$ and $\eps(n)n^{1/3} \to
\infty$. The following theorems, proved by Bollob\' as \cite{B1}
and \L uczak \cite{L} using different methods, describe the
behavior of the largest components when $p$ is outside the
``scaling-window''.

\begin{thm} {\bf [Subcritical phase]} \label{below} If $p(n)={1-\eps(n) \over n}$ then for any $\eta>0$ and integer $\ell>0$
we have $$ \prob \Big ( \Big | { |\C _\ell| \over
2\eps(n)^{-2}\log(n\eps(n)^3) } - 1
 \Big | > \eta \Big ) \to 0 \, ,$$
as $n \to \infty$.
\end{thm}
%

\begin{thm} {\bf [Supercritical phase]} \label{above} If $p(n)={1+\eps(n) \over n}$ then for any
$\eta >0$ we have
$$ \prob \Big ( \Big | { |\C _1| \over 2 n \eps(n) } - 1 \Big | >
\eta \Big ) \to 0 \, ,$$ and for any integer $\ell > 1$ we have
$$ \prob \Big ( \Big | { |\C _\ell| \over
2\eps^{-2}(n)\log(n\eps^3) } - 1 \Big | > \eta \Big ) \to 0 \, ,$$
as $n \to \infty$.
\end{thm}

The proofs of these theorems in \cite{B1} and \cite{L} are quite
involved and use the detailed asymptotics from \cite{Wr},
\cite{B1} and \cite{BCM} for the number of graphs on $k$ vertices
with $k+\ell$ edges. The proofs we present here are simple and
require no hard theorems. The main advantage, however, of these
proofs is their robustness. In a companion paper \cite{NP2} we use
similar methods to analyze critical percolation on a random
regular graphs. In this case, the enumerative methods employed in
\cite{B1} and \cite{L} are not available.

The phase transition in the Erd\H{o}s-R\'enyi random graphs
$G(n,p)$ happens when $p={c \over n}$. Namely, with high
probability, if $c>1$ then $|\C _1|$ is linear in $n$, and if
$c<1$ then $|\C_1|$ is logarithmic in $n$. When $c \sim 1$ the
situation is more delicate. In \cite{LPW}, \L uczak, Pittel and
Wierman prove that for $p={1 + \lambda n^{-1/3} \over n}$, the law
of $n^{-2/3} |\C_1|$ converges to a positive non-constant
distribution which in \cite{A} is identified as the longest
excursion length of some Brownian motion with variable drift. See
\cite{NP} for a recent account of the case $p={1 + \lambda
n^{-1/3} \over n}$ with simple proofs.

Thus, $|\C _1|$ is not concentrated and is roughly of size
$n^{2/3}$ if $p={1 + \lambda n^{-1/3} \over n}$. However, if
$\eps(n)$ a sequence such that $n^{1/3} \eps(n) \to \infty$ and
$p={1+ \eps(n) \over n}$ then as stated in Theorems \ref{below}
and \ref{above}, the size $|\C _1|$ of the largest component in
$G(n,p)$ is concentrated. In summary, $G(n,p)$ has a scaling
window of length $n^{-1/3}$ in which the percolation is
``critical'' in the sense that $|\C _1|$ is not concentrated.

\section{The exploration process}

We recall an exploration process, due to Karp and Martin-L\"of
(see \cite{K} and \cite{M}), in which vertices will be either {\em
active, explored} or {\em neutral}. After the completion of step
$t \in \{0,1,\ldots , n\}$ we will have precisely $t$ explored
vertices and the number of the active and neutral vertices is
denoted by $A_t$ and $N_t$ respectively.

Fix an ordering of the vertices $\{v_1, \ldots, v_n\}$. In step
$t=0$ of the process, we declare vertex $v_1$ {\em active} and all
other vertices {\em neutral}. Thus $A_0=1$ and $N_0 = n-1$. In
step $t \in \{1,\ldots,n\}$, if $A_{t-1}>0$ let $w_t$ be the first
active vertex; if $A_{t-1}=0$, let $w_t$ be the first neutral
vertex. Denote by $\eta_t$ the number of neutral neighbors of
$w_t$ in $G(n,p)$, and change the status of these vertices to {\em
active}. Then, set $w_t$ itself {\em explored}.

Denote by $\F _{t}$ the $\sigma$-algebra generated by $\{\eta_1,
\ldots , \eta_t\}$. Observe that given $\F_{t-1}$ the random
variable $\eta_t$ is distributed as Bin$(N_{t-1} - {\bf
1}_{\{A_{t-1}=0\}},p)$ and we have the recursions \be \label{nrec}
N_t = N_{t-1} - \eta_t - {\bf 1}_{\{A_{t-1}=0\}} \, , \qquad t
\leq n \, ,\ee and

\be \label{arec}
 A_t= \left \{
\begin{array}{ll}
A_{t-1} + \eta_t - 1, & A_{t-1} > 0 \\
\eta_t, & A_{t-1} = 0 \, , \qquad t
\leq n \, . \\
\end{array} \right .
\ee  As every vertex is either neutral, active or explored,

\be \label{nrec2} N_t = n - t - A_t \, , \qquad t \leq n \, . \ee

At each time $j\leq n$ in which $A_j = 0$, we have finished
exploring a connected component. Hence the random variable $Z_t$
defined by
$$ Z_t = \sum _{j=1}^{t-1} {\bf 1}_{\{A_j = 0\}} \, ,$$
counts the number of components completely explored by the process
before time $t$. Define the process $\{Y_t\}$ by $Y_0 = 1$ and
$$ Y_t = Y_{t-1} + \eta_t - 1 \, .$$ By (\ref{arec}) we have that
$Y_t = A_t - Z_t$, i.e. $Y_t$ counts the number of active vertices
at step $t$ minus the number of components completely explored
before step $t$.

At each step we marked as explored precisely one vertex. Hence,
the component of $v_1$ has size $\min \{t \geq 1 : A_t = 0\}$.
Moreover, let $t_1 < t_2 \ldots$ be the times at which
$A_{t_j}=0$; then $(t_1, t_2 - t_1, t_3 - t_2, \ldots)$ are the
sizes of the components. Observe that $Z_{t} = Z_{t_j}+1$ for all
$t \in \{t_j+1, \ldots, t_{j+1}\}$. Thus $Y_{t_{j+1}} = Y_{t_j}-1$
and if $t \in \{t_j +1, \ldots, t_{j+1}-1\}$ then $A_t
> 0$, and thus $Y_{t_{j+1}}< Y_t$. By induction we conclude that
$A_t=0$ if and only if $Y_t < Y_s$ for all $s<t$, i.e. $A_t=0$ if
and only if $\{Y_t\}$ has hit a new record minimum at time $t$. By
induction we also observe that $Y_{t_j} = -(j-1)$ and that for $t
\in \{t_j +1, \ldots t_{j+1}\}$ we have $Z_t = j$. Also, by our
previous discussion for $t \in \{t_j +1, \ldots t_{j+1}\}$ we have
$\min _{s \leq t-1} Y_t = Y_{t_j} = -(j-1)$, hence by induction we
deduce that $Z_t = - \min_{s \leq t-1} Y_t + 1$. Consequently, \be
\label{formula} A_t = Y_t - \min _{s\leq t-1} Y_s + 1 \, .\ee

\begin{lemma} \label{lowerbound} For all $p\leq {2 \over n}$ there
exists a constant $c>0$ such that for any integer $t>0$,
$$ \prob \Big (N_{t} \leq n - 5t \Big) \leq e^{-ct} \, .$$
\end{lemma}

\begin{proof} Let $\{\alpha_i\}_{i=1}^t$ be a sequence of i.i.d. random
variables distributed as Bin$(n,p)$. It is clear that we can
couple $\eta_i$ and $\alpha_i$ so $\eta_i \leq \alpha_i$ for all
$i$, and thus by (\ref{nrec}) \be \label{eq1} N_t \geq n-1 - t -
\sum _{i=1}^t \alpha _i \, .\ee The sum $\sum _{i=1}^t \alpha _i$
is distributed as Bin$(nt, p)$ and $p \leq {2 \over n}$ so by
Large Deviations (see \cite{AS} section A.14) we get that for some
fixed $c>0$
$$ \prob \Big ( \sum _{i=1}^t \alpha _i \geq 3t \Big ) \leq e^{-ct} \, ,$$
which together with (\ref{eq1}) concludes the proof.\end{proof}

\section{The subcritical phase}

Before beginning the proof of Theorem \ref{below} we require some
facts about processes with i.i.d. increments. Fix some small
$\eps>0$ and let $p={1-\eps \over m}$ for some integer $m>1$. Let
$\{\beta _j\}$ be a sequence of random variables distributed as
Bin$(m,p)$. Let $\{W_t\}_{t \geq 0}$ be a process defined by
$$ W_0 = 1, \qquad W_t = W_{t-1} + \beta_t - 1 \, .$$
Let $\tau$ be the hitting time of $0$,
$$ \tau = \min _t \{ W_t =  0\} \, .$$

By Wald's lemma we have that $\E \tau = \eps ^{-1}$. Further
information on the tail distribution of $\tau$ is given by the
following lemma.

\begin{lemma} \label{tail} There exists constant $C_1, C_2, c_1, c_2> 0$ such that for all $T>\eps^{-2}$ we have
$$ \prob (\tau \geq T ) \leq C_1\Big ( \eps^{-2} T^{-3/2} e^{- {(\eps^2 - c_1 \eps^3) T \over
2}}\Big ) \, ,$$ and
$$ \prob (\tau \geq T ) \geq c_1\Big ( \eps^{-2} T^{-3/2} e^{- {(\eps^2 + c_2 \eps^3) T \over
2}}\Big ) \, .$$ Furthermore,
$$ \E \tau^2 = O ( \eps^{-3} ) \, .$$
\end{lemma}

We will use the following proposition due to Spitzer (see
\cite{S}).

\begin{prop} \label{spitzer}
Let $a_0, \ldots, a_{k-1} \in {\mathbb Z}$ satisfy $\sum
_{i=0}^{k-1} a_i = -1$. Then there is precisely one $j \in \{0,
\ldots, k-1\}$ such that for all $r \in \{0, \ldots, k-2\}$
$$ \sum _{i=0}^{r} a_{(j+i) {\rm \  mod \ } k} \geq 0 \, .$$
\end{prop}

\noindent {\bf Proof of Lemma \ref{tail}.} By Proposition
\ref{spitzer}, $\prob (\tau = t) = {1 \over t} \prob (W_t = 0)$.
As $\sum _{j=1}^t \beta_j$ is distributed as a Bin$(mt,p)$ random
variable we have
$$ \prob (W_t = 0) = { mt \choose t-1}p^{t-1}(1-p)^{m-(t-1)} \, .$$
Replacing $t-1$ with $t$ in the above formula only changes it by a
multiplicative constant which is always between ${1/2}$ and $2$. A
straightforward computation using Stirling's approximation gives

\be \label{stirling} \prob (W_t = 0) = \Theta \Big \{ t^{-1/2}
(1-\eps)^t \Big (1+ {1 \over m-1} \Big )^{tm}\Big (1 - {1-\eps
\over m} \Big)^{t(m-1)} \Big \} \, .\ee Denote $ q = (1-\eps) \Big
( 1 + {1 \over m-1} \Big )^m \Big (1- {1-\eps \over m} \Big
)^{m-1} $, then
$$ \prob (\tau \geq T) = \sum _{t \geq T} \prob (\tau = t) = \sum _{t \geq
T}{1 \over t} \prob (W_t = 0) = \Theta \Big ( \sum _{t \geq T}
t^{-3/2} q^t \Big ) \, .$$ This sum can be bounded above by
$$ T^{-3/2} \sum _{t \geq T} q^t = T^{-3/2}{q^T \over 1-q} \, ,$$
and below by
$$ \sum _{t=T}^{2T} t^{-3/2} q^t \geq (2T)^{-3/2} {q^T(1-q^T)
\over 1-q} \, .$$ Observe that as $m \to \infty$ we have that $q$
tends to $(1-\eps)e^{\eps}$. By expanding $e^{\eps}$ we find that
$$ q = (1-\eps)(1+\eps+{\eps^2 \over 2}) + \Theta(\eps^3) = 1 -
{\eps^2 \over 2} + \Theta(\eps^3) \, .$$ Using this and the
previous bounds on $\prob (\tau \geq T)$ we get the first
assertion of the Lemma.

The second assertion follows from the following computation. By
(\ref{stirling}) we have that for some constant $C>0$
$$ \E \tau ^2 = \sum _{t\geq 1} t^2 \prob (\tau = t) = \sum_{t
\geq 1} t \prob (W_t = 0) \leq C\sum _{t \geq 1} \sqrt{t} q^t \,
.$$ Thus, by direct computation (or by \cite{F}, section XIII.5,
Theorem 5)
$$ \E \tau ^2 \leq O \Big ( {1 \over 1-q} \Big)^{3/2} =
O(\eps^{-3}) \, .$$ \qed \\


\noindent {\bf Proof of Theorem \ref{below}.} We begin with an
upper bound. Recall that component sizes are $t_{j+1} - t_j$ for
some $j>0$ where $t_j$ are record minima of the process $\{Y_t\}$.
For a vertex $v$ denote by $C(v)$ the connected component of
$G(n,p)$ which contains $v$. We first bound $\prob (|C(v_1)| > T)$
where
$$T =2(1+\eta)\eps^{-2}\log(n\eps^3) \, .$$
Recall that $|C(v_1)| = \min _t \{Y_t = 0\}$. Couple $\{Y_t\}$
with a process $\{W_t\}$ as in Lemma \ref{tail}, which has
increments distributed as Bin$(n,p)-1$ such that $Y_t \leq W_t$
for all $t$. Define $\tau$ as in Lemma \ref{tail}. As $p={1 - \eps
\over n}$ and $T> \eps^{-2}$, by Lemma \ref{tail} we have
$$ \prob (\tau > T) \leq C\eps (n\eps^3)^{-(1+\eta)} \log (n\eps^3)
^{-3/2}\, ,$$ for some fixed $C>0$. Our coupling implies that
$\prob ( |C(v_1)| > T ) \leq \prob ( \tau > T)$. Denote by $X$ the
number of vertices $v$ such that $|C(v)|>T$. If $|\C_1| > T$ then
$X>T$. Also, for any two vertices $v$ and $u$ by symmetry we have
that $|C(v)|$ and $|C(u)|$ are identically distributed. We
conclude that

\begin{eqnarray*} \prob ( |\C _1 | > T) &\leq& \prob (X > T) \leq
{\E X \over T} = { n \prob (|C(v_1)| > T) \over T} \\ &\leq& { C_1
n\eps (n\eps^3)^{-(1+\eta)(1-C_2\eps)} \log ^{-3/2} (n\eps^3)
\over 2(1+\eta)\eps^{-2}\log(n\eps^3)} \leq (n \eps
^3)^{{-\eta}(1-C_2\eps)+C_2\eps} \to 0   \, .
\end{eqnarray*}

We now turn to prove a lower bound. Write $$ T =
2(1-\eta)\eps^{-2}\log(n\eps^3) \, ,$$ and define the stopping
time
$$ \gamma = \min \{ t : N_t \leq n-{\eta \eps n \over 8} \} \, .$$ Recall
that $\{t_j\}$ are times in which $A_{t_j}=0$ and also $Y_{t_j}$
is a record minimum for $\{Y_t\}$.  For each integer $j$ let $\{
W^{(j)}_t\}$ be a process with increments distributed as
Bin$(n-{\eta \eps n \over 8}, p)$ where the starting point is
$W^{(j)}_0 = Y_{t_j} = -(j-1)$. Note that if $t_{j+1} < \gamma$
then we can couple $\{Y_t\}$ and $\{W^{(j)}_t\}$ such that $Y_{t_j
+ t} \geq W_t$ for all $t \in [t_j, t_{j+1}]$. Define the stopping
times $\{\tau_j\}$ by
$$ \tau_j = \min \{ t : W^{(j)}_t = -j \} \, .$$

\noindent Take
$$ N = \Big \lceil \eps^{-1} (n\eps^3)^{(1-{\eta \over 8})} \Big \rceil \,
.$$ We will prove that with high probability $t_N < \gamma$ and
that there exists $k_1 < k_2 < \ldots < k_\ell < N$ such that
$\tau _{k_i} > T$. Note that these two events imply that $|\C
_\ell| > T$. Indeed, by Lemma \ref{lowerbound} we have

\be \label{stp1} \prob \Big ( \gamma \leq {\eta \eps n \over 40}
\Big ) \leq e^{-c\eps n} \, .\ee By bounding the increments of
$\{Y_t\}$ above by variables distributed as Bin$(n,p)-1$ we learn
by Wald's Lemma (see \cite{D}) that $\E[t_{j+1} - t_j] \leq
\eps^{-1}$, hence $\E t_N \leq \eps^{-2} (n\eps^3)^{(1-{\eta \over
8})}$. We conclude that

\be \label{stp2} \prob ( t_N > {\eta \eps n \over 40} ) \leq { 40
\eps^{-2} (n\eps^3)^{(1-{\eta \over 8})} \over \eta \eps n } = {40
\over \eta}(n\eps^3)^{-{\eta \over 8}} \, ,\ee which goes to $0$
as $\eps n^{-1/3}$ tends to $\infty$. In Lemma \ref{tail} take
$m=n-{\eta \eps n \over 8}$ and note that $p={
(1-\eps)(1-{\eta\eps \over 8}) \over m} \geq { 1- (1+{\eta \over
8})\eps \over m}$, and so Lemma \ref{tail} gives that for any $j$
$$ \prob (\tau _j > T) \geq c_1 \eps (n \eps^3)^{-(1+{\eta \over 8})^2(1-\eta)(1+c_2\eps) }\log^{-3/2}(\eps^3 n) \geq \eps (n
\eps^3)^{-(1-{\eta \over 4})} \, .$$ Let $X$ be the number of $j
\leq N$ such that $\tau_j
> T$. Then we have
$$ \E X \geq N \eps (n \eps^3)^{-(1-{\eta \over 4})} \geq C
(n\eps^3)^{{\eta \over 8}} \to \infty \, ,$$ hence by Large
Deviations (see \cite{AS}, section A.14) for any fixed integer
$\ell>0$ we have
$$ \prob \Big ( X < \ell \Big ) \leq e^{-c (n\eps^3)^{{\eta \over
8}}} \, ,$$ for some fixed $c>0$. By our previous discussion, this
together with (\ref{stp1}) and (\ref{stp2}) gives
$$ \prob ( |\C _\ell| < T ) \leq O \Big ( {(n\eps^3)^{-{\eta \over
8}}\over \eta} \Big ) \, . $$
\qed \\

\section{The supercritical phase}

In this section we denote $\xi_t = \eta_t -1$. We first prove some
Lemmas.

\begin{lemma} \label{est1} If $p= {1+\eps \over n}$ then for all $t < 3\eps(n)n$
\be \label{est1.ineq1} \E A_t = O(\eps t + \sqrt{t}) \, ,\ee  and

\be \label{est1.ineq2} \E Z_t = O(\eps t + \sqrt{t}) \, .\ee
\end{lemma}

\begin{proof}
Write $T=3\eps n$. We will use (\ref{formula}). First observe that
as $\eta _t$ can always be bounded above by a Bin$(n,p)$ random
variable we can bound $\E \xi_t \leq \eps$ for all $t$ and hence
$\E Y_t \leq \eps t$. Denote by $\tau$ the stopping time $\tau =
\min \{ t : N_t \leq n - 15\eps n\}$. By definition of $\eta_t$ we
have
$$ \E [ \xi_t \mid \F_{t-1} ] = pN_{t-1} -p{\bf 1}_{\{A_{t-1}=0\}} -
1 \, .$$ As $\{N_t\}$ is a decreasing sequence, we deduce that as
long as $t < \tau$, we have $\E [ \xi_t \mid \F_{t-1} ] > -D\eps$
for $D>0$ large enough. Hence, the process $\{ D\eps j - Y_j \}
_{j=0}^{t \wedge \tau}$ is a submartingale for any $t$. By Doob's
maximal $L^2$ inequality we have

\be \label{doob} \E [ \max _{j \leq t\wedge\tau} (D\eps j - Y_j)^2
] \leq 4\E [ (D\eps (t\wedge \tau) -Y_{t \wedge \tau})^2] \, .\ee

For any $j < \tau$ the random variable $\eta_j$ can be
stochastically bounded from below by a Bin$(n-15\eps n,p)$ random
variable and above by a Bin$(n,p)$ random variable. Hence for any
$k < j < \tau$ we have
$$ \Big | \E [ \xi_j - D\eps \mid \F_{k} ] \Big | = O(\eps) \, ,$$
and so
$$ \E [ (\xi_j - D\eps)(\xi_k - D\eps) ] = O(\eps^2) \, .$$
We conclude that as long as $t < \tau$
$$ \E [ (D\eps t - Y_t)^2 ] \leq 2 \sum _{k< j}^t \E [(\xi_j -
D\eps)(\xi_k - D\eps)] + \sum _{j=1}^t \E[(\xi_j - D\eps)^2] =
O(\eps^2 t^2 + t) \, .$$ Lemma \ref{lowerbound} implies that for
$n$ large enough,

\be \label{LD} \prob \Big ( N_T \leq n - 15\eps n \Big ) \leq
e^{-3c\eps n} \leq {1 \over n^2} \, ,\ee and as $\{N_t\}$ is a
decreasing sequence we deduce that $\prob(\tau \leq T) \leq
n^{-2}$. Hence for any $t \leq T$
\begin{eqnarray*} \E [ (D\eps t -Y_{t})^2]&\leq& \E [ (D\eps (t\wedge \tau) -Y_{t \wedge
\tau})^2{\bf 1}_{\{t < \tau\}} ] + O(n^2)\prob (t \geq \tau) \\
&=& O(\eps^2 t^2 + t) \, .
\end{eqnarray*}
We deduce by (\ref{doob}) and Jensen inequality that for any $t
\leq T$
$$ \E [ \min _{j \leq t} (Y_j - D\eps j)] = O(\eps t + \sqrt{t})
\, ,$$ hence $\E [\min _{j \leq t} Y_j] = O(\eps t + \sqrt{t})$
and so by (\ref{formula}) we obtain (\ref{est1.ineq1}). Inequality
(\ref{est1.ineq2}) follows immediately from the relation $Z_t =
A_t - Y_t$.
\end{proof}

\begin{lemma} \label{est2} If $p={1 + \eps \over n}$ then for all $t < 3\eps(n)n$
\be \label{est2.ineq1} \E N_t = n(1-p)^t + O(\eps^2 n) \, ,\ee and

\be \label{est2.ineq2} \E \xi_t = \eps - {t \over n} + O(\eps^2)
\, .\ee
\end{lemma}

\begin{proof}
Observe that by (\ref{nrec}) we have that
$$ \E [ N_t \mid \F_{t-1} ] = (1-p)N_{t-1} - (1-p){\bf 1}_{\{A_{t-1}=0\}} \, .$$
By iterating this relation we get that $\E N_t = n(1-p)^t + O(\E
Z_t)$ which by Lemma \ref{est1} yields (\ref{est2.ineq1}) (observe
that for $t = 3\eps n$ we have $\eps t > \sqrt{t}$ by our
assumption on $\eps$). Since
$$E[\xi_t \mid \F_{t-1}] = pN_{t-1} - p {\bf 1}_{\{A_{t-1}=0\}}
-1\, ,$$ by taking expectations and using (\ref{est2.ineq1}) we
get

\begin{eqnarray*} \E \xi _t &=& (1+\eps)(1-{1+\eps \over n})^t -1 + O(\eps^2) \\
&=& (1+\eps)(1-(1+\eps)t/n) - 1 + O(\eps^2)  = \eps - {t \over n}
+ O(\eps^2) \, , \end{eqnarray*} where we used the fact that
$(1-x)^t = 1-tx+O(t^2 x^2)$.
\end{proof}

\noindent {\bf Proof of Theorem \ref{above}.}  Write $T=3\eps n$
and $\xi_j^* = \E [ \xi_j \mid \F_{j-1} ]$. The process
$$ M_t = Y_t - \sum _{j=1}^t \xi_j^* \, ,$$
is a martingale. By Doob's maximal $L^2$ inequality we have that
$$ \E ( \max _{t \leq T} M_t ^2) ) \leq 4\E M_T^2 \, . $$
As $M_t$ has orthogonal increments with bounded second moment we
conclude that $\E M_T^2 = O(T)$, hence, by Jensen's inequality we
have

\be \label{step1} \E \Big [\max _{t\leq T} \Big | Y_t - \sum
_{j=1}^t \xi_j^* \Big | \Big ] \leq O(\sqrt{T}) = O(\sqrt{\eps n})
\, . \ee
As $\xi_j^* = pN_{j-1} - p{\bf 1}_{\{A_{j-1}=0\}} -1$ by
(\ref{nrec2}) we have
$$ \E | \xi_j^*  - \E \xi_j | = p \E | A_{j-1} +
{\bf 1}_{\{A_{j-1}=0\}} - \E A_{j-1} - \E {\bf 1}_{\{A_{j-1}=0\}}
| \, .$$ By the triangle inequality and Lemma \ref{est1} we
conclude that for all $j\leq T$
$$ \E | \xi_j^*  - \E \xi_j | \leq p\cdot O(\eps j+\sqrt{j}) \, ,$$ and
hence for any $t \leq T$
$$ \E \Big [ \sum _{j \leq t}  | \xi_j^* - \E \xi_j
 | \Big ] \leq  p\cdot O(\eps t^2 + t^{3/2}) \leq O(\eps^3 n) \, .$$
By the triangle inequality we get

\be \label{step2} \E \Big [ \max _{t \leq T} \Big | \sum _{j=1}^t
( \xi_j^* - \E \xi_j) \Big | \Big ]\leq O(\eps^3 n)  \, . \ee
Using the triangle inequality, (\ref{step1}), (\ref{step2}) and
Markov inequality gives

\be \label{keystep} \prob \Big ( \max _{t\leq T} \Big | Y_t - \sum
_{j=1}^t \E \xi_j \Big | \geq \delta \eps^2 n \Big ) \leq
\delta^{-1} ( O(\eps) + O ( (\eps^3 n)^{-1/2} )) \longrightarrow 0
\, . \ee


\noindent Lemma \ref{est2} implies that for any $b>0$

\begin{eqnarray} \label{key3}
\sum _{j=1}^{b\eps n} \E \xi_j = \sum _{j=1}^{b\eps n} \Big (\eps
- {t \over n} + O(\eps^2)\Big ) = (b-{b^2 \over 2})\eps^2 n +
O(\eps^3 n) \, .
%
\end{eqnarray}
By (\ref{keystep}) and (\ref{key3}) we deduce that for $\eta>0$
small enough, with probability tending to $1$, the process $Y_t$
is strictly positive at times $[\eta \eps n, (2-\eta) \eps n]$ and
hence
$$ \prob \Big ( |\C _1| > 2(1-\eta)\eps n \Big ) \geq 1 - O \Big (
\eta^{-1}(\eps + (\eps^3 n)^{-1/2}) \Big ) \, .$$ We also deduce
by (\ref{keystep}) and (\ref{key3}) that at time $t= (2+\eta)\eps
n$ we have $Y_t \leq -{\eta^2 \over 3} \eps ^2 n$ and at all times
$t<\eta \eps n$ we have that $Y_t > - {\eta^2 \over 3} \eps^2 n$
with probability tending to $1$. As component sizes are excursion
lengths of $Y_t$ above its past minima, we conclude that by time
$2(1+\eta)\eps n$ we have explored completely at least one
component of size at least $2(1-\eta)\eps n$. As $N_t \leq n-t$
for $t>2(1-\eta)\eps n$ and $\eta < 1/4$ we have $\E[ \eta_t -1
\mid \F_{t-1} ] \leq -{\eps \over 2}$. By optional stopping we
immediately get that if $t_j > 2(1-\eta)\eps n$ then $\E [t_{j+1}
- t_j] \leq 2\eps^{-1}$. Thus if $\C$ is a component which we
began discovering after time $2(1-\eta)\eps n$ we have
$$ \prob ( |\C| \geq \eps n ) \leq {2 \over \eps^2 n} \, .$$
Denote by $X(\eta)$ the number of vertices of which $|C(v)| > \eps
n$ which we began discovering after time $2(1-\eta)\eps n$, then
we learn that $\E X(\eta) \leq {2 \over \eps ^2}$. Denote by $\C_1
(\eta)$ the largest component which we began discovering after
time $2(1-\eta)\eps n$. Clearly if $|\C_1(\eta)| > \eps n$ then
$X(\eta) > \eps n$, thus by Markov inequality
$$ \prob (|\C_1(\eta)| > \eps n) \leq {2 \over \eps^3 n} \to 0
\, .$$ Thus we have proved that there exists a unique component of
size between $2(1-\eta)\eps n$ and $2(1+\eta)\eps n$. Condition on
this event and consider the graph remained on the complement of
this component. This graph has $m$ vertices where
$$ | m - (n-2\eps n) | < 2\eta \eps n \, ,$$
and as $p={1+\eps \over n}$ we have that
$$ \Big | p - \Big ( {1 - {\eps \over m}} \Big) \Big | \leq {2 \eta \eps + O(\eps^2) \over m } \, .$$
This graph is distributed as $G(m,p)$ restricted to the event that
it does not contain a component of size between $2(1-\eta)\eps n$
and $2(1+\eta)\eps n$. By Theorem \ref{below} we know that the
event that there exists such a component has probability $o(1)$.
Thus for any collection of graphs ${\cal B}$ on $m$ vertices which
does not contain such a component, the probability of ${\cal B}$
in the remaining graph is $(1+o(1)) \prob_{m,p} ({\cal B})$ where
$\prob _{m,p}$ is the usual $G(m,p)$ probability measure. Thus, we
conclude by Theorem \ref{below} that for any integer $\ell > 1$
and $\eta'>0$
$$ \prob \Big ( \Big | { |\C _\ell| \over
2\eps^{-2}(n)\log(n\eps^3) } - 1
 \Big | > \eta' \Big ) \to 0 \, ,$$ concluding the
proof of the theorem. \qed \newline

{\bf Remark.} With a little more effort it is possible to show for
the supercritical case, that in the exploration process for any
fixed $\ell$, the $\ell$-th largest component is discovered {\em
after} the largest component is discovered.

\section*{Acknowledgments}

The first author would like to thank Microsoft Research, in which
parts of this research were conducted, for their kind hospitality.

The work was done partially while the authors were visiting the
Institute for Mathematical Sciences, National University of
Singapore in 2006. The visit was supported by the Institute.

\bigskip \noindent
{\bf Asaf Nachmias}: \texttt{asafnach(at)math.berkeley.edu} \\
Department of Mathematics\\
UC Berkeley\\
Berkeley, CA 94720, USA.

\bigskip \noindent
{\bf Yuval Peres}: \texttt{peres(at)stat.berkeley.edu} \\
Microsoft Research\\
One Microsoft way,\\
Redmond, WA 98052-6399, USA.

\end{document}